\documentclass[11pt]{amsart}
\usepackage{xcolor,mathtools}

\include{package}
\definecolor{grn}{rgb}{0,0.6,0}
\definecolor{mrn}{rgb}{0.3,0,0}
\definecolor{blue}{rgb}{0,0,0.7}
\definecolor{Mygray}{rgb}{0.75,0.75,0.75}
\definecolor{auburn}{rgb}{0.43, 0.21, 0.1}
\definecolor{britishracinggreen}{rgb}{0.0, 0.26, 0.15}
\definecolor{taupe}{rgb}{0.28, 0.24, 0.2}
\usepackage{amsmath,amssymb}
\newtheorem{theorem}{Theorem}

\newtheorem{proposition}{Proposition}[section]

\newtheorem{lemma}{Lemma}

\newtheorem{remark}{Remark}

\newtheorem*{ack}{Acknowledgements}

\usepackage{latexsym,amssymb}
\usepackage{amssymb}
\usepackage{graphicx}
\usepackage{ textcomp }
\usepackage{tikz}
\usetikzlibrary{shapes.geometric, arrows}
\tikzstyle{startstop} = [rectangle, rounded corners, minimum width=3cm, minimum height=1cm,text centered, draw=black, fill=white!30]
\tikzstyle{io} = [trapezium, trapezium left angle=70, trapezium right angle=110, minimum width=3cm, minimum height=1cm, text centered, draw=black, fill=white!30]
\tikzstyle{process} = [rectangle, minimum width=2cm, minimum height=1cm, text centered, draw=black, fill=white!30]
\tikzstyle{decision} = [rectangle, minimum width=1cm, minimum height=1cm, text centered, draw=black, fill=white!30]
\tikzstyle{arrow} = [thick,->,>=stealth]
\usepackage[left=2cm,right=2cm,top=2cm,bottom=2cm]{geometry}

\allowdisplaybreaks

\begin{document}
\baselineskip=14.5pt
\title[Bi-quadratic P\'{o}lya fields and large P\'{o}lya groups]{Some bi-quadratic P\'{o}lya fields and large P\'{o}lya groups of compositum of simplest cubic and quintic fields}

\author{Md. Imdadul Islam, Debopam Chakraborty and Jaitra Chattopadhyay}
\address[Md. Imdadul Islam and Debopam Chakraborty]{Department of Mathematics, BITS-Pilani, Hyderabad campus, Hyderabad, INDIA}
\address[Jaitra Chattopadhyay]{Department of Mathematics, Siksha Bhavana, Visva-Bharati, Santiniketan - 731235, West Bengal, India}

\email[Md. Imdadul Islam]{p20200059@hyderabad.bits-pilani.ac.in}

\email[Debopam Chakraborty]{debopam@hyderabad.bits-pilani.ac.in}

\email[Jaitra Chattopadhyay]{jaitra.chattopadhyay@visva-bharati.ac.in; chat.jaitra@gmail.com}

\begin{abstract}
The P\'{o}lya group $Po(K)$ of an algebraic number field $K$ is the subgroup of the ideal class group $Cl_{K}$ generated by the ideal classes of the products of prime ideals of the same norm. If $Po(K)$ is trivial, then the number field $K$ is said to be a P\'{o}lya field. In this article, we furnish three families $\mathbb{Q}(\sqrt{p},\sqrt{qrs})$, $\mathbb{Q}(\sqrt{2p},\sqrt{qrs})$ and $\mathbb{Q}(\sqrt{2p},\sqrt{2qrs})$ of bi-quadratic P\'{o}lya fields $K$ involving prime numbers $p,q,r$ and $s$ that satisfy certain quadratic residue conditions. It is worthwhile to note that in each of the fields, exactly five primes ramify in $K/\mathbb{Q}$ and this is the maximum possible number of ramified primes in a P\'{o}lya field over $\mathbb{Q}$. Towards the end of the paper, we discuss about large P\'{o}lya groups of the compositums of Shank's cubic fields and Lehmer's quintic fields and prove that there are infinitely many such fields with index $1$.
\end{abstract}

\renewcommand{\thefootnote}{}

\footnote{2020 \emph{Mathematics Subject Classification}: Primary 11R29, Secondary 11R11.}

\footnote{\emph{Key words and phrases}: P\'{o}lya fields, P\'{o}lya groups, Galois cohomology.}

\footnote{\emph{We confirm that all the data are included in the article.}}

\renewcommand{\thefootnote}{\arabic{footnote}}
\setcounter{footnote}{0}

\maketitle

\section{introduction}

Let $K$ be an algebraic number field with ring of integers $\mathcal{O}_{K}$, discriminant $d_{K}$, ideal class group $Cl_{K}$ and class number $h_{K}$. The {\it Ostrowski ideals} (cf. \cite{ost}) of $K$ are defined by \begin{equation*}
\displaystyle\Pi_{p^{f}}(K) = \displaystyle\prod_{\substack {\mathfrak{p} \in {\rm{Max}}(\mathcal{O}_{K})\\ N_{K/\mathbb{Q}} (\mathfrak{m}) = p^{f}}} \mathfrak{m}. 
\end{equation*}

If all the Ostrowski ideals of for any prime power are principal ideals, then $K$ is said to be a {\it P\'{o}lya field} (cf. \cite{pol}, \cite{zan}). As $Cl_{K}$ measures the failure of $\mathcal{O}_{K}$ from being a principal ideal domain, the P\'{o}lya group, which is the subgroup of $Cl_{K}$ generated by the ideal classes of all the Ostrowski ideals of $K$ (cf. \cite{cah1}, \textsection II.4), measures the failure of $K$ from being a P\'{o}lya field. Thus $K$ is a P\'{o}lya field if and only if $Po(K)$ is trivial. Due to the works of Brumer and Rosen in \cite{bru}, it is known that if $K/\mathbb{Q}$ is a Galois extension then $Po(K)$ is freely generated by the ideal classes of the ambiguous ideals. 

\smallskip

Many analogous questions regarding $Cl_{K}$ turn out to be much easier to address for $Po(K)$. In some cases, it is answered affirmatively. For example, Golod and Shafarevich proved in \cite{golod} that there exist number fields that cannot be embedded inside a number field of class number $1$. But Leriche \cite{ler-em} proved that any number field can indeed be embedded inside a P\'{o}lya field. In fact, she proved that the Hilbert class field of any number field is a P\'{o}lya field.

\smallskip

Another elusive question in algebraic number theory is the Gauss class number one problem which asserts that there are infinitely many real quadratic fields of class number one. Till date it is unknown whether there exists infinitely many number fields of class number one. However, a complete classification of quadratic, cubic and sextic P\'{o}lya fields is known due to the works of Hilbert \cite[Theorem 106]{hil} and Leriche \cite{ler}. Leriche also characterized bi-quadratic P\'{o}lya fields that are compositum of two quadratic P\'{o}lya fields in \cite[Theorem 5.4]{ler}.

\smallskip

There have been quite a few studies (cf. \cite{cha1}, \cite{cha2}, \cite{hei1}, \cite{hei2}, \cite{isl}, \cite{isl-rnt}, \cite{tou}) on the P\'{o}lya groups of totally real bi-quadratic fields. In this paper, we focus on some families of totally real bi-quadratic P\'{o}lya fields. More precisely, we prove the following theorems. 

\begin{theorem}\label{mainthm-5.1}
Let $p,q,r,s$ be prime numbers with $p \equiv q \equiv r \equiv 3 \pmod 4, qrs \equiv 1 \pmod 4$. Assume that $\left(\frac{q}{s}\right) = \left(\frac{r}{s}\right) = 1$ does not hold true. Then the bi-quadratic field $\mathbb{Q}(\sqrt{p}, \sqrt{qrs})$ is a P\'{o}lya field if $\left(\frac{p}{s}\right) \neq \left(\frac{2}{s}\right) = -1$, and $\left(\frac{qr}{p}\right) = \left(\frac{pr}{q}\right) = -1$. Moreover, the condition $\left(\frac{pr}{q}\right) = -1$ is redundant if $\left(\frac{q}{s}\right) = \left(\frac{r}{s}\right) = -1$.    
\end{theorem}

\begin{theorem}\label{mainthm-5.2.1}
      Let $p,q,r,s$ be primes such that $p \equiv q \equiv 3 \pmod 4, r \equiv s \equiv 1 \pmod 4$. Then $\mathbb{Q}(\sqrt{p}, \sqrt{qrs})$ is a P\'{o}lya field in the following cases.
    \subitem $(i)$ Given $\left(\frac{2}{r}\right) = \left(\frac{2}{s}\right) = -1$, if any of the following conditions hold. \\
    C1: $\left(\frac{q}{r}\right) = \left(\frac{q}{s}\right) = \left(\frac{rs}{p}\right) -1$. \\
    C2: $\left(\frac{rs}{q}\right) = \left(\frac{r}{p}\right) = \left(\frac{s}{p}\right) =  -1$.     
    \subitem $(ii)$ Given $\left(\frac{2}{r}\right) = -\left(\frac{2}{s}\right) = 1$, if any of  the following conditions hold. \\
    C3: $\left(\frac{r}{s}\right) = \left(\frac{q}{s}\right) = \left(\frac{r}{p}\right) = -1 = -\left(\frac{r}{q}\right)$. \\
    C4: $\left(\frac{r}{q}\right) = \left(\frac{r}{s}\right) = \left(\frac{r}{p}\right) = -1 = -\left(\frac{s}{p}\right)$.\\ 
    C5: $\left(\frac{r}{q}\right) = \left(\frac{r}{p}\right) = \left(\frac{pq}{s}\right) = -1 = -\left(\frac{r}{s}\right)$.
\end{theorem}
\begin{remark}\label{rmk-5.2.1}
    We note that in the statement of Theorem \ref{mainthm-5.2.1}, $\left(\frac{2}{r}\right) = -\left(\frac{2}{s}\right) = -1$ is also a possibility. But due to $r \equiv s \pmod 4$, a simple renaming of those two primes will yield the same result. 
\end{remark}

\begin{theorem}\label{mainthm-5.2.2}
      Let $p,q,r,s$ be primes such that $p \equiv q \equiv 3 \pmod 4, r \equiv s \equiv 3 \pmod 4$. Then $\mathbb{Q}(\sqrt{p}, \sqrt{qrs})$ is a P\'{o}lya field in the following cases.
    \subitem $(i)$ $\left(\frac{2}{r}\right) \neq \left(\frac{2}{s}\right)$, $\left(\frac{q}{r}\right) = \left(\frac{q}{s}\right) = \left(\frac{p}{q}\right) = 1 = -\left(\frac{p}{r}\right) = -\left(\frac{p}{s}\right)$.
    \subitem $(ii)$ $\left(\frac{2}{r}\right) = \left(\frac{2}{s}\right) \neq \left(\frac{2}{q}\right)$, $\left(\frac{q}{r}\right) = \left(\frac{r}{s}\right) = \left(\frac{s}{q}\right)$, $\left(\frac{pr}{s}\right) = \left(\frac{ps}{r}\right) = -1$.
\end{theorem}
\begin{remark}
    Before we write down the proof of the theorem, we note that the result captures all possible cases except the case $q \equiv r \equiv s \pmod 8$, in a manner similar to Remark \ref{rmk-5.2.1}.
\end{remark}
\begin{theorem}\label{mainthm-5.3}
      Let $p,q,r,s$ be primes such that $p \equiv q \equiv r \equiv 3 \pmod 4, qrs \equiv 1 \pmod 4$ and $\left(\frac{q}{s}\right) = \left(\frac{r}{s}\right) = 1$ does not hold true. Then the bi-quadratic field $\mathbb{Q}(\sqrt{2p}, \sqrt{qrs})$ is a P\'{o}lya field in the following cases.
    \subitem $(i)$ Given $\left(\frac{q}{s}\right) = \left(\frac{r}{s}\right) = -1$, if any of the following conditions hold. \\
    C1: $\left(\frac{2p}{s}\right) = \left(\frac{qr}{p}\right) = -\left(\frac{2p}{q}\right) = -1$, $q \not \equiv r \pmod 8$. \\
    C2: $\left(\frac{s}{p}\right) = -\left(\frac{qr}{p}\right) = 1$, $p \equiv q \equiv r \pmod 8$, $s \equiv 5 \pmod 8$.\\
    C3: $\left(\frac{s}{p}\right) = \left(\frac{qr}{p}\right) = -1$, $p \not \equiv q \equiv r \pmod 8$, $s \equiv 1 \pmod 8$.
 
    \subitem $(ii)$ Given $\left(\frac{q}{s}\right) \neq \left(\frac{r}{s}\right) = 1$, if any of  the following conditions hold. \\
    C4: Assuming $\left(\frac{r}{p}\right) = \left(\frac{p}{q}\right) = \left(\frac{q}{r}\right) = -1$, either $\left(\frac{2}{p}\right) = \left(\frac{2}{r}\right) = \left(\frac{2}{s}\right) = -1$, $\left(\frac{s}{p}\right) = 1$ or $\left(\frac{2}{p}\right) = \left(\frac{2}{q}\right) = - \left(\frac{2}{r}\right) = \left(\frac{2}{s}\right) = 1$, $\left(\frac{s}{p}\right) = -1$.\\
    C5: Assuming $\left(\frac{r}{p}\right) = \left(\frac{p}{q}\right) = \left(\frac{r}{q}\right) = -1$, either $\left(\frac{2}{p}\right) = -\left(\frac{2}{r}\right) = -\left(\frac{2}{s}\right) = 1$, $\left(\frac{s}{p}\right) = 1$ or $\left(\frac{2}{p}\right) \neq \left(\frac{2}{q}\right)$, $\left(\frac{2}{r}\right) = -1$, $\left(\frac{2}{s}\right) = 1$, $\left(\frac{s}{p}\right) = -1$.\\
    C6: Assuming $\left(\frac{p}{r}\right) = \left(\frac{q}{p}\right) = \left(\frac{r}{q}\right) = -1$, either $\left(\frac{2}{p}\right) = \left(\frac{2}{r}\right) = -\left(\frac{2}{s}\right) = 1$, $\left(\frac{s}{p}\right) = 1$ or $\left(\frac{2}{p}\right) = \left(\frac{2}{q}\right) = - \left(\frac{2}{r}\right) = -\left(\frac{2}{s}\right) = -1$, $\left(\frac{s}{p}\right) = -1$.\\
    C7: Assuming $\left(\frac{p}{r}\right) = \left(\frac{q}{p}\right) = \left(\frac{q}{r}\right) = -1$, either $\left(\frac{2}{p}\right) = -\left(\frac{2}{r}\right) = \left(\frac{2}{s}\right) = -1$, $\left(\frac{s}{p}\right) = 1$ or $\left(\frac{2}{p}\right) \neq \left(\frac{2}{q}\right)$, $\left(\frac{2}{r}\right) = \left(\frac{2}{s}\right) = 1$, $\left(\frac{s}{p}\right) = -1$.

\end{theorem}

\begin{remark}
    We note that in the statement of Theorem \ref{mainthm-5.3}, $\left(\frac{q}{s}\right) \neq \left(\frac{r}{s}\right) = -1$ is also a possibility. But due to $q \equiv r \equiv 3 \pmod 4$, a simple renaming of those two primes will yield the same result. 
\end{remark}
\begin{theorem}\label{mainthm-5.4.1}
      Let $p,q,r,s$ be primes such that $p \equiv q \equiv 3 \pmod 4, r \equiv s \equiv 1 \pmod 4$. Then the bi-quadratic field $\mathbb{Q}(\sqrt{2p}, \sqrt{2qrs})$ is a P\'{o}lya field if $$\left(\frac{2}{s}\right) = \left(\frac{q}{r}\right) = \left(\frac{s}{r}\right) = \left(\frac{p}{r}\right) = -1 = -\left(\frac{q}{s}\right) = -\left(\frac{p}{s}\right).$$
\end{theorem}
\begin{theorem}\label{mainthm-5.4.2}
      Let $p,q,r,s$ be primes such that $p \equiv q \equiv 3 \pmod 4, r \equiv s \equiv 3 \pmod 4$. Then the bi-quadratic field $\mathbb{Q}(\sqrt{2p}, \sqrt{2qrs})$ is a P\'{o}lya field if $$\left(\frac{2}{r}\right) = \left(\frac{2}{s}\right) \neq \left(\frac{2}{q}\right), \left(\frac{q}{r}\right) = \left(\frac{r}{s}\right) = \left(\frac{s}{q}\right), \left(\frac{pr}{s}\right) = \left(\frac{ps}{r}\right) = -1.$$
\end{theorem}

\section{preliminaries}

If $K/\mathbb{Q}$ is a finite Galois extension with Galois group $G$. Then the group of units $\mathcal{O}_{K}^{*}$ is naturally endowed with a $G$-module structure through the action $(\sigma,\alpha) \mapsto \sigma(\alpha)$. This makes it possible to obtain a relation between the cohomology group $H^{1}(G,\mathcal{O}_{K}^{*})$ and the ramified primes in $K/\mathbb{Q}$. The following proposition of Zantema provides one such relation.
\begin{proposition}\label{prop-zan}\cite[Page 163]{zan}
Let $G$ be the Galois group of the finite Galois extension $K/\mathbb{Q}$ and let $e_{1},\ldots,e_{s}$ be the ramification indices of the ramified primes in $K/\mathbb{Q}$. Then we have an exact sequence 
    \begin{equation}\label{exact-equn}
0 \to H^{1}(G,\mathcal{O}_{K}^{*}) \to \displaystyle\bigoplus_{i = 1}^{s}\mathbb{Z}/e_{i}\mathbb{Z} \to Po(K) \to 0
\end{equation}
of finite abelian groups. 
\end{proposition}

It is often easier to understand the group $H^{1}(G,\mathcal{O}_{K}^{*})$ via some of its subgroups. The following results of Setzer and Zantema help us achieve that. The first one furnishes a necessary and sufficient condition of $H^{1}(G,\mathcal{O}_{K}^{*})$ being equal to its $2$-torsion subgroup $H^{1}(G,\mathcal{O}_{K}^{*}[2])$. The precise statement is as follows.
\begin{lemma}\label{lem-zan1}\cite[Theorem 4]{set}
Let $K$ be a totally real bi-quadratic field and let $K_{1}, K_{2}$ and $K_{3}$ be the quadratic subfields of $K$. Then $[H^{1}(G,\mathcal{O}_{K}^{*}) : H^{1}(G,\mathcal{O}_{K}^{*})[2]] \leq 2$. Equality holds if and only if the rational prime $2$ is totally ramified in $K/\mathbb{Q}$ and there exists $\gamma_{i} \in K_{i}$ for each $i = 1,2,3$ with $N(\gamma_{1}) = N(\gamma_{2}) = N(\gamma_{3}) = \pm 2$. 
\end{lemma}
The next theorem due to Zantema facilitates us to compute the group $H^{1}(G,\mathcal{O}_{K}^{*}[2])$ by identifying it with a subgroup of $\mathbb{Q}^{*}/(\mathbb{Q}^{*})^{2}$. For an integer $t \neq 0$, let $[t]$ denote its image in the group $\mathbb{Q}^{*}/(\mathbb{Q}^{*})^{2}$. Also, for non-zero integers $t_{1},\ldots,t_{r}$, let $\langle [t_{1}],\ldots,[t_{r}] \rangle$ stand for the subgroup of $\mathbb{Q}^{*}/(\mathbb{Q}^{*})^{2}$ generated by their canonical images.
\begin{lemma}\label{lem-zan2}\cite[Lemma 4.3]{zan}
Let $K$ be a totally real bi-quadratic field and let $K_{1}, K_{2}$ and $K_{3}$ be the quadratic subfields of $K$. For $i=1,2,3$, let $\Delta_{i}$ be the square-free part of the discriminant of $K_{i}$ and let $u_{i} = x_{i}+y_{i}\sqrt{\Delta_{i}}$ be a fundamental unit of $\mathcal{O}_{K_{i}}$ with $x_{i} > 0$. Then $H^{1}(G, \mathcal{O}_{K}^{*})[2]$ is isomorphic to the subgroup $\langle [\Delta_{1}], [\Delta_{2}], [\Delta_{3}], [a_{1}], [a_{2}], [a_{3}] \rangle \subseteq \mathbb{Q}^{*}/ (\mathbb{Q}^{*})^{2}$ where $a_{i} = N(u_{i}+1)$ if $N(u_{i}) = 1$, and $1$ otherwise.
\end{lemma}
\begin{lemma}\label{lem5.1}
Let $ K_{2} = \mathbb{Q}(\sqrt{qrs})$ where $q, r$ and $s$ are prime numbers with $q \equiv r \equiv 3 \pmod 4$, $s \equiv 1 \pmod 4$ and such that $qrs \equiv 1 \pmod 4$. Let $u_{2} = x_{2} + y_{2} \sqrt{qrs}$ be a fundamental unit of $K_{2}$ and let $a_{2} = N(u_{2} + 1)$. Then the square-free representative of $[a_{2}]$ is a proper, non-trivial divisor of $qrs$. More precisely, we have the following.
\begin{enumerate}
\item If $\left(\frac{q}{r}\right) = \left(\frac{q}{s}\right) = 1$, then $[a_{2}] = [q], [s]$ or $[qs]$. Moreover, if $\left(\frac{r}{s}\right) = -1$, then $[a_{2}] = [q]$. 
\item If $\left(\frac{q}{r}\right) = \left(\frac{q}{s}\right) = -\left(\frac{r}{s}\right)  = -1$, then $[a_{2}] = [r]$. 
\item If $\left(\frac{q}{s}\right) = \left(\frac{r}{s}\right) =  -1$, then $[a_{2}] = [qr]$. 
\item If $\left(\frac{q}{r}\right) = - \left(\frac{q}{s}\right) = \left(\frac{r}{s}\right) = 1$, then $[a_{2}] = [qs]$ . 
\item If $\left(\frac{q}{r}\right) = - \left(\frac{q}{s}\right) = -1$, then $[a_{2}] = [r], [s]$ or $[rs]$. Moreover, if $\left(\frac{r}{s}\right) = -1$, then $[a_{2}] = [rs]$.
\end{enumerate}
\end{lemma}

The proof of this is exactly the same as Lemma 2.4 of \cite{isl-rnt}. We give a tabular description for better understanding of the lemma below. 

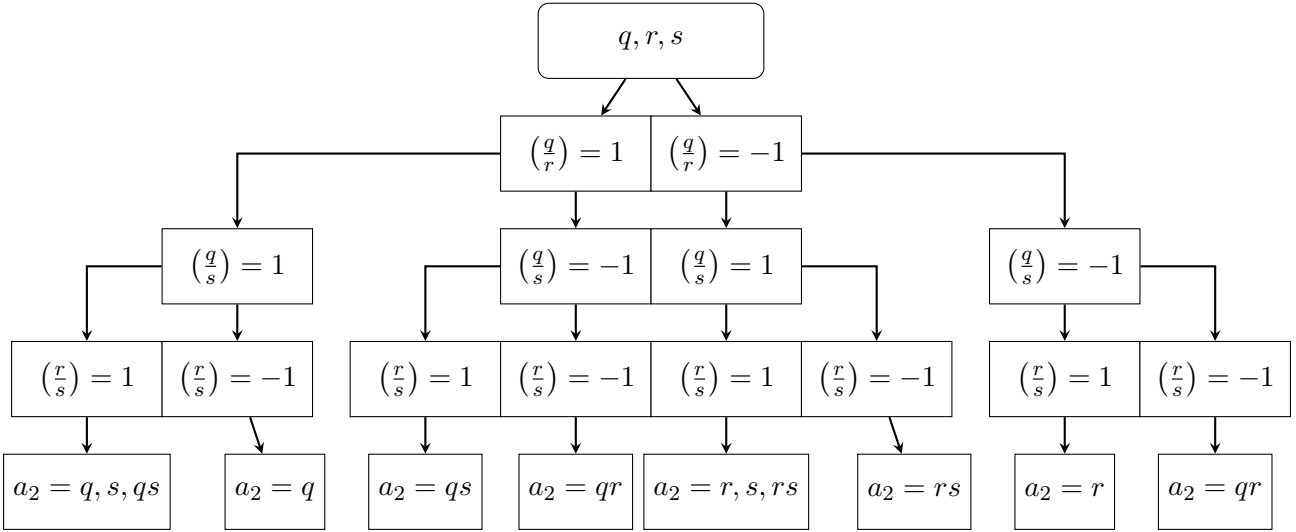
\begin{center}
\begin{tikzpicture}[node distance=1.5cm]
\node (start) [startstop] {$q, r, s$};
\node (pro2a) [process, below of=start, xshift=-1cm] {$\left(\frac{q}{r}\right) = 1$};
\node (pro2b) [process, below of=start, xshift=1cm] {$\left(\frac{q}{r}\right) = -1$};
\node (pro3a) [process, below of=pro2a, xshift=-4.5cm] {$\left(\frac{q}{s}\right) = 1$};
\node (pro3b) [process, below of=pro2a, xshift=0cm] {$\left(\frac{q}{s}\right) = -1$};
\node (pro3c) [process, below of=pro2b, xshift=0cm] {$\left(\frac{q}{s}\right) = 1$};
\node (pro3d) [process, below of=pro2b,xshift=4.5cm] {$\left(\frac{q}{s}\right) = -1$};
\node (pro4a) [process, below of=pro3a, xshift=-2cm] {$\left(\frac{r}{s}\right) = 1$};
\node (pro4b) [process, below of=pro3a, xshift=0cm] {$\left(\frac{r}{s}\right) = -1$};
\node (pro4c) [process, below of=pro3b, xshift=-2cm] {$\left(\frac{r}{s}\right) = 1$};
\node (pro4d) [process, below of=pro3b, xshift=0cm] {$\left(\frac{r}{s}\right) = -1$};
\node (pro4e) [process, below of=pro3c, xshift=-0cm] {$\left(\frac{r}{s}\right) = 1$};
\node (pro4f) [process, below of=pro3c, xshift=2cm] {$\left(\frac{r}{s}\right) = -1$};
\node (pro4g) [process, below of=pro3d, xshift=0cm] {$\left(\frac{r}{s}\right) = 1$};
\node (pro4h) [process, below of=pro3d, xshift=2cm] {$\left(\frac{r}{s}\right) = -1$};
\node (dec1) [decision, below of=pro4a, xshift=0cm] {$a_{2} = q, s, qs$};
\node (dec2) [decision, below of=pro4b, xshift=0.5cm] {$a_{2} = q$};
\node (dec3) [decision, below of=pro4c, xshift=0cm] {$a_{2} = qs$};
\node (dec4) [decision, below of=pro4d, xshift=0cm] {$a_{2} = qr$};
\node (dec5) [decision, below of=pro4e, xshift=0cm] {$a_{2} = r, s, rs$};
\node (dec6) [decision, below of=pro4f, xshift=0.5cm] {$a_{2} = rs$};
\node (dec7) [decision, below of=pro4g, xshift=0cm] {$a_{2} = r$};
\node (dec8) [decision, below of=pro4h, xshift=0cm] {$a_{2} = qr$};
\draw [arrow] (start) -- (pro2a);
\draw [arrow] (start) -- (pro2b);
\draw [arrow] (pro2a) -| (pro3a);
\draw [arrow] (pro2a) -- (pro3b);
\draw [arrow] (pro2b) -- (pro3c);
\draw [arrow] (pro2b) -| (pro3d);
\draw [arrow] (pro3a) -| (pro4a);
\draw [arrow] (pro3a) -- (pro4b);
\draw [arrow] (pro3b) -| (pro4c);
\draw [arrow] (pro3b) -- (pro4d);
\draw [arrow] (pro3c) -- (pro4e);
\draw [arrow] (pro3c) -| (pro4f);
\draw [arrow] (pro3d) -- (pro4g);
\draw [arrow] (pro3d) -| (pro4h);
\draw [arrow] (pro4a) -- (dec1);
\draw [arrow] (pro4b) -- (dec2);
\draw [arrow] (pro4c) -- (dec3);
\draw [arrow] (pro4d) -- (dec4);
\draw [arrow] (pro4e) -- (dec5);
\draw [arrow] (pro4f) -- (dec6);
\draw [arrow] (pro4g) -- (dec7);
\draw [arrow] (pro4h) -- (dec8);
\end{tikzpicture}    
\end{center}

\begin{remark}
    We note that $qrs \equiv 5 \pmod 8$ might lead to a fundamental unit with $2$ in the denominator. A quick calculation can yield that the outcomes for $[a_{2}]$ still remain unchanged in the above context. Hence, the choice of $qrs \equiv 1 \pmod 4$ was mentioned without breaking into two cases modulo $8$.
\end{remark}

\begin{remark}\label{rem5}
Prior to formulating the proofs, we also note the following result due to Maarefparvar. If $u$ denotes the fundamental unit of a quadratic field $\mathbb{Q}(\sqrt{d})$, the result helps us compute two fewer possibilities of $N(u+1)$.

\begin{proposition}[\cite{mar1}, Proposition 2.1]\label{prop-mar}
Let $k= \mathbb{Q}(\sqrt{d})$ be a real quadratic field with fundamental unit $u > 1$. Assume that $N_{k/ \mathbb{Q}}(u) = 1$ and let $n_{k}$ denote the square-free part of $N_{k/ \mathbb{Q}}(u+1)$. Then $n_{k} \not \in \{1,d\}$, $n_{k}$ divides the discriminant of $k$, $n_{k}$ is the norm of an integer in $k$, and $n_{k} \cdot u$ is a square in $k$.
\end{proposition}
    
\end{remark}

\section{Proof of Theorem \ref{mainthm-5.1}}

We first note that $[\Delta_{1}] = p, [\Delta_{2}] = qrs, [\Delta_{3}] = pqrs$. Also, by Remark \ref{rem5} and subsequent Proposition \ref{prop-mar}, we have $[a_{1}] = 2 \mbox{ or } 2p$. Let $u_{2} = x_{2} + y_{2} \sqrt{qrs}$ be the fundamental unit of $K_{2} = \mathbb{Q}(\sqrt{qrs})$. Then $[a_{2}] = N(u_{2}+1)$. Given $q \equiv r \equiv 3 \pmod 4$, $N(u_{2}) = x_{2}^{2} - qrs y_{2}^{2} = 1$ implies that $N(u_{2}+1) = 2(x_{2}+1)$ and $x_{2}$ is always odd. This implies $[a_{2}] = \frac{x_{2}+1}{2}$. A detailed description of $[a_{2}]$ is already mentioned in Lemma \ref{lem5.1}. From Lemma \ref{lem5.1}, we also note that, if $\left(\frac{q}{s}\right) = \left(\frac{r}{s}\right) = -1$ then $[a_{2}] = [qr]$ and this implies $\langle{\Delta_{1}, \Delta_{2}, \Delta_{3}, [a_{1}], [a_{2}]\rangle} = \langle{2,p, qr, s\rangle}.$ Let $u_{3} = x_{3} + y_{3} \sqrt{pqrs}$ be the fundamental unit of $K_{3} = \mathbb{Q}(\sqrt{pqrs})$. Then $[a_{3}] = N(u_{3}+1)$. Given $q \equiv 3 \pmod 4$, $N(u_{3}) = x_{3}^{2} - pqrs y_{3}^{2} = 1$ implies that $N(u_{3}+1) = 2(x_{3}+1)$. This implies $[a_{3}] = 2(x_{3}+1)$ if $x_{3}$ is even, and $[a_{3}] = \frac{x_{3}+1}{2}$ if $x_{3}$ is odd.

\smallskip

Moreover, we also have $(x_{3}+1)(x_{3}-1) = pqrsy_{2}^{2}$ which implies $x_{3}+1 = t_{1}a^{2}, x_{3}-1 = t_{2}b^{2}$ such that $t_{1}t_{2} = pqrs$ and $ab = y_{3}$ when $x_{3}$ is even (respectively, $\frac{x_{3}+1}{2} = t_{1}a^{2}, \frac{x_{3}-1}{2} = t_{2}b^{2}$ such that $t_{1}t_{2} = pqrs$ and $ab = \frac{y_{3}}{2}$, when $x_{3}$ is odd). 

\smallskip

For $x_{3}$ even, that is, $x_{3}+1$ odd, we note from above that $t_{1}a^{2} - t_{2}b^{2} = 2, \text{ } t_{1}t_{2} = pqrs.$
We note that $t_{1} \not \in \{1, pqrs\}$ as $\left(\frac{2}{s}\right) = -1$. As $\left(\frac{2}{s}\right) \neq \left(\frac{qr}{s}\right)$, one can note that $t_{1} \not \in \{ps, qr\}$. Moreover, $\left(\frac{2}{s}\right) \neq \left(\frac{p}{s}\right)$ implies $t_{1} \not \in \{p, s, qrs, pqr\}$. Now, if $x_{3}$ is odd, we solve $t_{1}a^{2} - t_{2}b^{2} = 1, \text{ } t_{1}t_{2} = pqrs$. Given $t_{1} \equiv 3 \equiv -t_{2} \pmod 4$ implies $a^{2} + b^{2} \equiv -1 \pmod 4$, a contradiction, we note that $t_{1} \not \in \{p, ps, pqr\}$. As $\left(\frac{qr}{p}\right) = -1$, one can immediately observe that $t_{1} \neq qr$. Moreover, as $\left(\frac{p}{s}\right) = \left(\frac{s}{p}\right) = 1$, we note that $t_{1} \neq qrs$ either. Also, $\left(\frac{q}{s}\right) = \left(\frac{s}{q}\right) = -1$ implies $t_{1} \neq s$ here. 

\smallskip

Together, this shows that if $\left(\frac{q}{s}\right) = \left(\frac{r}{s}\right) = -1$, then $\langle{\Delta_{1}, \Delta_{2}, \Delta_{3}, [a_{1}], [a_{2}], [a_{3}]\rangle} = \langle{2,p, q, r, s\rangle}$. That is, $K$ is P\'{o}lya. Also $\left(\frac{q}{s}\right) \neq \left(\frac{r}{s}\right)$ implies $[a_{2}] \in \{[r], [qs]\}$ or $\{[q], [rs]\}$, depending on whether $\left(\frac{r}{s}\right) = 1$ or $-1$, respectively. This, in turn, implies $\langle{\Delta_{1}, \Delta_{2}, \Delta_{3}, [a_{1}], [a_{2}]\rangle} = \langle{2,p, qs, r\rangle}$ if $\left(\frac{q}{s}\right) \neq \left(\frac{r}{s}\right) = 1$, and $\langle{\Delta_{1}, \Delta_{2}, \Delta_{3}, [a_{1}], [a_{2}]\rangle} = \langle{2,p, rs, q\rangle}$ if $\left(\frac{q}{s}\right) \neq \left(\frac{r}{s}\right) = -1$. The proof for these two cases with the additional condition of $\left(\frac{pr}{q}\right) = -1$ follows a very similar approach. By using Proposition \ref{prop-zan} and Lemma \ref{lem-zan2}, it follows that $\mathbb{Q}(\sqrt{p},\sqrt{qrs})$ is a P\'{o}lya field. 

\section{Proof of Theorem \ref{mainthm-5.2.1}}

In this case, we only prove the first sub-case in detail, as the other cases follow a very similar approach. We note that $[\Delta_{1}] = p, [\Delta_{2}] = qrs, [\Delta_{3}] = pqrs, [a_{1}] = 2 \text{ or } 2p$. Let $u_{2} = x_{2} + y_{2} \sqrt{qrs}$ be the fundamental unit of $K_{2} = \mathbb{Q}(\sqrt{qrs})$. Then $[a_{2}] = N(u_{2}+1)$. Given $q \equiv 3 \pmod 4$, $N(u_{2}) = x_{2}^{2} - qrs y_{2}^{2} = 1$ implies that $N(u_{2}+1) = 2(x_{2}+1)$. This implies $[a_{2}] = \frac{x_{2}+1}{2}$ if $x_{2}$ is odd, and $2(x_{2}+1)$ if $x_{2}$ is even. As $(x_{2}+1) (x_{2}-1) = qrs y_{2}^{2}$, we note that $[a_{2}] = t_{1}$ when $x_{2}$ is odd (respectively, $2t_{1}$ when $x_{2}$ is even) where $t_{1}a^{2} - t_{2}b^{2} = 1, t_{1}t_{2} = qrs, ab = \frac{y_{2}}{2}$ when $x_{2}$ is odd (respectively, $t_{1}a^{2} - t_{2}b^{2} = 2, t_{1}t_{2} = qrs, ab = y_{2}$ when $x_{2}$ is even).  

\smallskip

Now, for $x_{2}$ is odd, we first note that for $t_{1} \in \{q, qs, qr\}$, the above equation looks like $a^{2}+b^{2} \equiv -1 \pmod 4$, hence have no solution. Also, for $t_{1} \in \{r, s, sr\}$, the condition $\left(\frac{q}{r}\right) = \left(\frac{q}{s}\right) = -1$ ensures there are no solution for those cases either. 

\smallskip

For $x_{2}$ even, we now solve $$t_{1}a^{2} - t_{2}b^{2} = 2, t_{1}t_{2} = qrs, ab = y_{2}.$$ Given $\left(\frac{2}{r}\right) = -1$, one can immediately observe that $t_{1} \not \in \{1, qrs\}$. For $t_{1} = r$, the equation is $r a^{2} - qs b^{2} = 2$. This shows $\left(\frac{2}{q}\right) = \left(\frac{r}{q}\right)$, $\left(\frac{2}{s}\right) = \left(\frac{r}{s}\right)$. We also note that $\left(\frac{2}{q}\right) = \left(\frac{r}{q}\right) = -1$. Hence, $1 = \left(\frac{2}{r}\right) \left(\frac{2}{s}\right) = \left(\frac{qs}{r}\right) \left(\frac{r}{s}\right) = \left(\frac{q}{r}\right) \left(\frac{s}{r}\right)^{2} = \left(\frac{q}{r}\right) = -1$, a contradiction. Hence, $t_{1} \neq r$. A very similar argument also shows that $t_{1} \not \in \{s, qr, qs\}$. Hence, $[a_{2}] \in \{2q, 2rs\}$, which then implies $\langle{\Delta_{1}, \Delta_{2}, \Delta_{3}, [a_{1}], [a_{2}]\rangle} = \langle{2,p, q, rs\rangle}$.

\smallskip

Let $u_{3} = x_{3} + y_{3} \sqrt{pqrs}$ be the fundamental unit of $K_{3} = \mathbb{Q}(\sqrt{pqrs})$. Then $[a_{3}] = N(u_{3}+1)$. Given $q \equiv 3 \pmod 4$, $N(u_{3}) = x_{3}^{2} - pqrs y_{3}^{2} = 1 \implies N(u_{3}+1) = 2(x_{3}+1)$. Now $pqrs \equiv 1 \pmod 4$ implies that $x_{3}$ must always be odd. This, in turn, shows that $[a_{3}] = \frac{x_{3}+1}{2}$. Denoting $\frac{x_{3}+1}{2}$ and $\frac{x_{3}-1}{2}$ by $t_{1}$ and $t_{2}$ respectively, we now try to solve $$t_{1}a^{2} - t_{2}b^{2} = 1, t_{1}t_{2} = pqrs, ab = \frac{y_{3}}{2}.$$
We first note that $\left(\frac{rs}{p}\right) = -1$ implies $\left(\frac{r}{p}\right) \neq \left(\frac{s}{p}\right)$. This immediately implies $t_{1} \not \in \{p, rs, qrs\}$. For $t_{1} = pq$, the equation is $pqa^{2} - rsb^{2} = 1$ which implies $\left(\frac{pq}{r}\right) = \left(\frac{pq}{s}\right)$. As $\left(\frac{q}{r}\right) = \left(\frac{q}{s}\right)$, this, in turn, implies $\left(\frac{p}{r}\right) = \left(\frac{p}{s}\right)$, a contradiction again. Now, $\left(\frac{q}{r}\right) = \left(\frac{q}{s}\right) = -1$ implies $t_{1} \not \in \{q, r, s, prs\}$. Together, this implies $$\langle{\Delta_{1}, \Delta_{2}, \Delta_{3}, [a_{1}], [a_{2}], [a_{3}]\rangle} = \langle{2,p, q, r, s\rangle}.$$ Consequently, by Proposition \ref{prop-zan} and Lemma \ref{lem-zan2} it follows that $\mathbb{Q}(\sqrt{p},\sqrt{qrs})$ is a P\'{o}lya field. This concludes the proof of C1. 

\smallskip

For C2, we note that $\left(\frac{r}{q}\right) \neq \left(\frac{s}{q}\right) = 1$ implies $[a_{2}] \in \{s, 2s, 2qr\}$, while $\left(\frac{r}{p}\right) = \left(\frac{s}{p}\right) = -1$ enforces that $[a_{3}] \not \in \langle{2, p, qr, s\rangle}$. Similarly, $\left(\frac{r}{q}\right) \neq \left(\frac{s}{q}\right) = -1$ implies $[a_{2}] \in \{r, 2r, 2qs\}$, while $\left(\frac{r}{p}\right) = \left(\frac{s}{p}\right) = -1$ enforces that $[a_{3}] \not \in \langle{2, p, qs, r\rangle}$. In both cases, it proves $\mathbb{Q}(\sqrt{p},\sqrt{qrs})$ is a P\'{o}lya field. 

\smallskip

For C3, under the given conditions, we note that, $[a_{2}] \in \{2q, 2rs\}$, while $[a_{3}] \not \in \langle{2, p, q, rs\rangle}$ implying $\mathbb{Q}(\sqrt{p},\sqrt{qrs})$ is a P\'{o}lya field. Similarly, for C4, $[a_{2}] \in \{2r, 2qs\}$, and $[a_{3}] \not \in \langle{2, p, r, qs\rangle}$, hence $\mathbb{Q}(\sqrt{p},\sqrt{qrs})$ is a P\'{o}lya field. And finally for C5, $[a_{2}] \in \{s, 2s, 2qr\}$, $[a_{3}] \not \in \langle{2, p, s, qr\rangle}$ together again implies $\mathbb{Q}(\sqrt{p},\sqrt{qrs})$ is a P\'{o}lya field.

\section{Proof of Theorem \ref{mainthm-5.2.2}}

Here, we also only prove the first sub-case, as the other cases can be proven via very similar approaches. We first note that $[\Delta_{1}] = p, [\Delta_{2}] = qrs, [\Delta_{3}] = pqrs, [a_{1}] = 2 \text{ or } 2p$. Let $u_{2} = x_{2} + y_{2} \sqrt{qrs}$ is the fundamental unit of $K_{2} = \mathbb{Q}(\sqrt{qrs})$. Then $[a_{2}] = N(u_{2}+1)$. Given $q \equiv 3 \pmod 4$, $N(u_{2}) = x_{2}^{2} - qrs y_{2}^{2} = 1$ implies that $N(u_{2}+1) = 2(x_{2}+1)$. This implies $[a_{2}] = \frac{x_{2}+1}{2}$ if $x_{2}$ is odd, and $2(x_{2}+1)$ if $x_{2}$ is even.

\smallskip

As $(x_{2}+1) (x_{2}-1) = qrs y_{2}^{2}$, we note that $[a_{2}] = t_{1}$ when $x_{2}$ is odd (resp. $2t_{1}$ when $x_{2}$ is even) where $t_{1}a^{2} - t_{2}b^{2} = 1, t_{1}t_{2} = qrs, ab = \frac{y_{2}}{2}$ when $x_{2}$ is odd (resp. $t_{1}a^{2} - t_{2}b^{2} = 2, t_{1}t_{2} = qrs, ab = y_{2}$ when $x_{2}$ is even).  \\
\noindent Now for $x_{2}$ is odd, we first note that for $t_{1} \in \{q, r, s\}$, the above equation looks like $a^{2}+b^{2} \equiv -1 \pmod 4$, hence have no solution. Also, for $t_{1} = rs$, the condition $1 = \left(\frac{q}{r}\right) = -\left(\frac{-q}{r}\right)$ ensures there are no solution for that case either. This implies $t_{1} \in \{qr, qs\}$ here.\\
\noindent For $x_{2}$ even, we now solve $$t_{1}a^{2} - t_{2}b^{2} = 2, t_{1}t_{2} = qrs, ab = y_{2}.$$ Given $\left(\frac{2}{r}\right) \neq \left(\frac{2}{s}\right)$, one can immediately observe that $t_{1} \not \in \{1, qrs\}$. For $t_{1} = r$, the equation is $r a^{2} - qs b^{2} = 2$. This shows $\left(\frac{2}{r}\right) = \left(\frac{-qs}{r}\right) = -\left(\frac{q}{r}\right)\left(\frac{s}{r}\right) = -\left(\frac{s}{r}\right) = \left(\frac{r}{s}\right)$, $\left(\frac{2}{s}\right) = \left(\frac{r}{s}\right)$. As $\left(\frac{2}{r}\right) \neq \left(\frac{2}{s}\right)$, we now get, $-1 = \left(\frac{2}{r}\right) \left(\frac{2}{s}\right) = \left(\frac{r}{s}\right)^{2} = 1$, a contradiction. Hence, $t_{1} \neq r$. A very similar argument also shows that $t_{1} \not \in \{q, s, qr, qs, rs\}$.\\
\noindent Hence, $[a_{2}] \in \{qr, qs\}$, and we get $$\langle{\Delta_{1}, \Delta_{2}, \Delta_{3}, [a_{1}], [a_{2}]\rangle} = \text{ either } \langle{2,p, qr, s\rangle} \text{ or } \langle{2,p, qs, r\rangle}.$$ 
Now, let $u_{3} = x_{3} + y_{3} \sqrt{pqrs}$ be the fundamental unit of $K_{3} = \mathbb{Q}(\sqrt{pqrs})$. Then $[a_{3}] = N(u_{3}+1)$. Given $q \equiv 3 \pmod 4$, $N(u_{3}) = x_{3}^{2} - pqrs y_{3}^{2} = 1 \implies N(u_{3}+1) = 2(x_{3}+1)$. Now $pqrs \equiv 1 \pmod 4$ implies that $x_{3}$ must always be odd. This, in turn, shows that $[a_{3}] = \frac{x_{3}+1}{2}$. Denoting $\frac{x_{3}+1}{2}$ and $\frac{x_{3}-1}{2}$ by $t_{1}$ and $t_{2}$ respectively, we now try to solve $$t_{1}a^{2} - t_{2}b^{2} = 1, t_{1}t_{2} = pqrs, ab = \frac{y_{3}}{2}.$$
We first note that $\left(\frac{p}{r}\right) = -\left(\frac{r}{p}\right) = -1$ implies $t_{1} \not \in \{p, pqs\}$. Similarly, $\left(\frac{p}{q}\right) = 1$ implies $t_{1} \neq qrs$, and $\left(\frac{q}{r}\right) = -\left(\frac{r}{q}\right) = 1$ implies $t_{1} \neq r$, while $\left(\frac{q}{s}\right) = -\left(\frac{s}{q}\right) = 1$ contradicts that $t_{1} = s$. Now $\left(\frac{p}{q}\right) = 1 = \left(\frac{q}{r}\right)$ implies $t_{1} \neq pr$, while $\left(\frac{p}{q}\right) = 1 = -\left(\frac{p}{s}\right)$ implies $t_{1} \neq qs$. Now $\left(\frac{p}{s}\right) = -\left(\frac{s}{p}\right) = -1$ implies $t_{1} \neq pqr$. Moreover, $\left(\frac{p}{q}\right) = 1 = \left(\frac{q}{s}\right)$ implies $t_{1} \neq ps$, while $\left(\frac{p}{q}\right) = 1 = -\left(\frac{p}{r}\right)$ implies $t_{1} \neq qr$. Together, this implies $[a_{3}] \in \{q, pq, rs, prs\}$. Hence, $$\langle{\Delta_{1}, \Delta_{2}, \Delta_{3}, [a_{1}], [a_{2}], [a_{3}]\rangle} = \langle{2,p, q, r, s\rangle}.$$ Hence by Proposition \ref{prop-zan} and Lemma \ref{lem-zan2} it follows that $\mathbb{Q}(\sqrt{p},\sqrt{qrs})$ is a P\'{o}lya field.

\section{Proof of Theorem \ref{mainthm-5.3}}

Similar to Theorem \ref{mainthm-5.2.1}, here also we only prove the first sub-case, as the other cases can be proven via very similar approaches. We first note that $[\Delta_{1}] = 2p, [\Delta_{2}] = qrs, [\Delta_{3}] = 2pqrs, [a_{1}] = 2 \text{ or } p$. Let $u_{2} = x_{2} + y_{2} \sqrt{qrs}$ is the fundamental unit of $K_{2} = \mathbb{Q}(\sqrt{qrs})$. Then $[a_{2}] = N(u_{2}+1)$. Given $q \equiv 3 \pmod 4$, and $qrs \equiv 1 \pmod 4$, $N(u_{2}) = x_{2}^{2} - qrs y_{2}^{2} = 1 \implies N(u_{2}+1) = 2(x_{2}+1)$ and $x_{2}$ is always odd. This implies $[a_{2}] = \frac{x_{2}+1}{2}$.\\

\noindent From Lemma \ref{lem5.1}, we note that, if $\left(\frac{q}{s}\right) = \left(\frac{r}{s}\right) = -1$ then $[a_{2}] = [qr]$ which implies $$\langle{\Delta_{1}, \Delta_{2}, \Delta_{3}, [a_{1}], [a_{2}]\rangle} = \langle{2,p, qr, s\rangle}.$$ 
Let $u_{3} = x_{3} + y_{3} \sqrt{2pqrs}$ be the fundamental unit of $K_{3} = \mathbb{Q}(\sqrt{2pqrs})$. Then $[a_{3}] = N(u_{3}+1)$. Given $q \equiv 3 \pmod 4$, $N(u_{3}) = x_{3}^{2} - 2pqrs y_{3}^{2} = 1 \implies N(u_{3}+1) = 2(x_{3}+1)$ and $x_{3}$ is odd. This implies $[a_{3}] = \frac{x_{3}+1}{2}$.\\
Moreover, we also have $(x_{3}+1)(x_{3}-1) = 2pqrsy_{2}^{2}$ which implies $\frac{x_{3}+1}{2} = t_{1}a^{2}, \frac{x_{3}-1}{2} = t_{2}b^{2}$, i.e., $$ t_{1}a^{2} - t_{2}b^{2} = 1, \text{ where } t_{1}t_{2} = pqrs, \text{ and } ab = \frac{y_{3}}{2}.$$ 
We first focus on the case $\left(\frac{q}{s}\right) = \left(\frac{r}{s}\right) = -1$. We first note that from the given conditions, the following are true in all three cases. $$\left(\frac{qr}{p}\right) = \left(\frac{2p}{s}\right) = -1, \text{ and } \left(\frac{2}{l}\right) = -1 \text{ for at least one } l \in \{p,q,r,s\}.$$
The above conditions along with $\left(\frac{q}{s}\right) = \left(\frac{r}{s}\right) = -1$ implies $t_{1} \not \in \{2, s, 2p, qr, qrs, 2pqr, pqrs\}$. For $t_{1} = p$, the equation $pa^{2} - 2qrs b^{2} = 1$ which then implies $\left(\frac{p}{q}\right) = \left(\frac{p}{r}\right) \implies \left(\frac{qr}{p}\right) = 1$, a contradiction. Hence $t_{1} \neq p$. Similarly, as $\left(\frac{q}{s}\right) = \left(\frac{r}{s}\right)$, $2qr a^{2} - ps b^{2} = 1 \implies \left(\frac{ps}{q}\right) = \left(\frac{ps}{r}\right)$ which again implies $\left(\frac{p}{q}\right) = \left(\frac{p}{r}\right)$ as above, a contradiction. Hence, $t_{1} \neq 2qr$. A very similar approach to the previous two cases also yields in $t_{1} \not \in \{ps, 2qrs\}$. \\
\noindent For $t_{1} = 2ps$, we note that the equation is $2ps a^{2} - qr b^{2} = 1$. Given $\left(\frac{q}{s}\right) = \left(\frac{r}{s}\right)$, we get that $\left(\frac{2ps}{q}\right) = \left(\frac{2ps}{r}\right) \implies \left(\frac{2p}{q}\right) = \left(\frac{2p}{r}\right)$ which, in turn, implies $\left(\frac{2}{q}\right) \neq \left(\frac{2}{r}\right)$ as $\left(\frac{qr}{p}\right) = -1$. This can only happen in the case $q \not \equiv r \pmod 8$, but that implies $\left(\frac{2p}{q}\right) = 1 \implies \left(\frac{2ps}{q}\right) = -1$, a contradiction. Hence, $t_{1} \neq 2ps$.\\
For equation $pqr a^{2} - 2s b^{2} = 1$, i.e., $t_{1} = pqr$, we first note that $q \not \equiv r \pmod 8$ implies $\left(\frac{-2s}{q}\right) \neq \left(\frac{-2s}{r}\right)$, a contradiction. A similar argument also holds true when $p \not \equiv q \pmod 8$. Now when $p \equiv q \equiv r \pmod 8$, $1 = \left(\frac{pqr}{s}\right) \implies \left(\frac{p}{s}\right) = \left(\frac{s}{p}\right) = 1$ as $\left(\frac{qr}{s}\right) = 1$. But then $1 = \left(\frac{-2s}{p}\right) \implies \left(\frac{2}{p}\right) = -1 \implies \left(\frac{2}{q}\right) = -1$ as $p \equiv q \pmod 8$. This implies $1 = \left(\frac{-2s}{q}\right) = \left(\frac{-1}{q}\right) \cdot \left(\frac{2}{q}\right) \cdot \left(\frac{s}{q}\right) = -1$, a contradiction again. Hence $t_{1} \neq pqr$. \\
for $t_{1} = 2s$, the equation is $2s a^{2} - pqr b^{2} = 1$. Given $\left(\frac{qr}{s}\right) = 1$, this immediately implies that $\left(\frac{p}{s}\right) = \left(\frac{s}{p}\right) = 1$, a contradiction, if we consider the third sub-case. Now under the assumption that $\left(\frac{s}{p}\right) = 1$, we note that then $1 = \left(\frac{2s}{p}\right) \implies \left(\frac{2}{p}\right) = 1$. Now $1 = \left(\frac{2s}{q}\right)$ implies $\left(\frac{2}{q}\right) = -1$ as $\left(\frac{s}{q}\right) = -1$. This implies $p \not \equiv q  \pmod 8$, hence a contradiction if the second sub-case is considered. Finally, under the assumptions of the first sub-case, we note that $\left(\frac{2s}{q}\right) = \left(\frac{2s}{r}\right) \implies \left(\frac{2}{q}\right) = \left(\frac{2}{r}\right)$ as $\left(\frac{q}{s}\right) = \left(\frac{r}{s}\right)$. But, this also leads to a contradiction as $q \not \equiv r \pmod 8$ in this case. Hence, $t_{1} \neq 2s$ either. 

\smallskip

Now, after eliminating all possible cases in $\langle{2, p, qr, s \rangle}$, we can conclude that $$\langle{\Delta_{1}, \Delta_{2}, \Delta_{3}, [a_{1}], [a_{2}], [a_{3}]\rangle} = \langle{2, p, q, r, s\rangle}.$$ Therefore, by Proposition \ref{prop-zan} and Lemma \ref{lem-zan2} it follows that $\mathbb{Q}(\sqrt{2p},\sqrt{qrs})$ is a P\'{o}lya field.

\section{Proof of Theorem \ref{mainthm-5.4.1}} 

We first note that $[\Delta_{1}] = 2p, [\Delta_{2}] = 2qrs, [\Delta_{3}] = pqrs, [a_{1}] = 2 \text{ or } p$. Let $u_{2} = x_{2} + y_{2} \sqrt{2qrs}$ is the fundamental unit of $K_{2} = \mathbb{Q}(\sqrt{2qrs})$. Then $[a_{2}] = N(u_{2}+1)$. Given $q \equiv 3 \pmod 4$, $N(u_{2}) = x_{2}^{2} - 2qrs y_{2}^{2} = 1 \implies N(u_{2}+1) = 2(x_{2}+1)$. This implies $[a_{2}] = \frac{x_{2}+1}{2}$ and $x_{2}$ is always odd.\\
\noindent As $(x_{2}+1) (x_{2}-1) = 2qrs y_{2}^{2}$, we note that $[a_{2}] = t_{1}$, where $t_{1}a^{2} - t_{2}b^{2} = 1, t_{1}t_{2} = 2qrs, ab = \frac{y_{2}}{2}$. Now, we first note that $\left(\frac{2}{s}\right) = -1$ implies $t_{1} \not \in \{2, qrs\}$. Similarly, $\left(\frac{q}{r}\right) = -1$ implies $t_{1} \not \in \{q, r , 2rs\}$ whereas $\left(\frac{s}{r}\right) = -1$ implies $t_{1} \not \in \{s, 2qr, 2qs\}$. The equation for $t_{1} = 2q$ is $2q a^{2} - rs b^{2} = 1$. But $\left(\frac{2}{s}\right) = -1$, $\left(\frac{q}{s}\right) = 1$ implies $\left(\frac{2q}{s}\right) = -1$, a contradiction. Hence, $t_{1} \neq 2q$. Similarly, $\left(\frac{rs}{q}\right) = -1$ implies $t_{1} \neq rs$, and $\left(\frac{qr}{s}\right) = -1$ leads to $t_{1} \not \in \{2s, qr\}$. Hence, $[a_{2}] \in \{2r, qs\}$. Consequently,
$$\langle{\Delta_{1}, \Delta_{2}, \Delta_{3}, [a_{1}], [a_{2}]\rangle} = \langle{2,p, qs, r\rangle}.$$ 
Let $u_{3} = x_{3} + y_{3} \sqrt{pqrs}$ be the fundamental unit of $K_{3} = \mathbb{Q}(\sqrt{pqrs})$. Then $[a_{3}] = N(u_{3}+1)$. Given $q \equiv 3 \pmod 4$, $N(u_{3}) = x_{3}^{2} - pqrs y_{3}^{2} = 1 \implies N(u_{3}+1) = 2(x_{3}+1)$. Denoting $\frac{x_{3}+1}{2}$ and $\frac{x_{3}-1}{2}$ by $t_{1}$ and $t_{2}$ respectively, we now try to solve $$t_{1}a^{2} - t_{2}b^{2} = 1, t_{1}t_{2} = pqrs, ab = \frac{y_{3}}{2}.$$
\noindent Now, we first note that $\left(\frac{p}{r}\right) = -1$ implies $t_{1} \not \in \{p, qrs\}$. Similarly, $\left(\frac{q}{r}\right) = -1$ implies $t_{1} \neq q$ whereas $\left(\frac{r}{s}\right) = -1$ implies $t_{1} \not \in \{r, s, pqs\}$. Now, $\left(\frac{pr}{s}\right) = -1$ implies $t_{1} \not \in \{pr, qs\}$. Hence, $[a_{3}] \in \{pq, ps, rs, qr, pqr, prs\}$. Together, this implies $$\langle{\Delta_{1}, \Delta_{2}, \Delta_{3}, [a_{1}], [a_{2}], [a_{3}]\rangle} = \langle{2,p, q, r, s\rangle}.$$ Thus by Proposition \ref{prop-zan} and Lemma \ref{lem-zan2} it follows that $\mathbb{Q}(\sqrt{2p},\sqrt{2qrs})$ is a P\'{o}lya field.

\section{Proof of Theorem \ref{mainthm-5.4.2}}

We first note that $[\Delta_{1}] = 2p, [\Delta_{2}] = 2qrs, [\Delta_{3}] = pqrs, [a_{1}] = 2 \text{ or } p$. Let $u_{2} = x_{2} + y_{2} \sqrt{2qrs}$ is the fundamental unit of $K_{2} = \mathbb{Q}(\sqrt{2qrs})$. Then $[a_{2}] = N(u_{2}+1)$. Given $q \equiv 3 \pmod 4$, $N(u_{2}) = x_{2}^{2} - 2qrs y_{2}^{2} = 1 \implies N(u_{2}+1) = 2(x_{2}+1)$. This implies $[a_{2}] = \frac{x_{2}+1}{2}$ and also, $x_{2}$ is always odd.\\
\noindent As $(x_{2}+1) (x_{2}-1) = 2qrs y_{2}^{2}$, we note that $[a_{2}] = t_{1}$, where $t_{1}a^{2} - t_{2}b^{2} = 1, t_{1}t_{2} = 2qrs, ab = \frac{y_{2}}{2}$. Now, we first note that $\left(\frac{2}{s}\right) \neq \left(\frac{2}{q}\right) $ implies $t_{1} \not \in \{2, qrs\}$. Similarly, $\left(\frac{q}{r}\right) = \left(\frac{r}{s}\right)$ implies $t_{1} \not \in \{r , qs, 2qs\}$ whereas $\left(\frac{q}{r}\right) = \left(\frac{s}{q}\right)$ implies $t_{1} \not \in \{q, rs, 2rs\}$. Moreover, $\left(\frac{r}{s}\right) = \left(\frac{s}{q}\right)$ implies $t_{1} \not \in \{s, qr, 2qr\}$. The equation for $t_{1} = 2q$ is $2q a^{2} - rs b^{2} = 1$. But $\left(\frac{2}{s}\right) = \left(\frac{2}{s}\right)$ while $\left(\frac{q}{r}\right) \neq \left(\frac{q}{s}\right)$, implying $\left(\frac{2q}{r}\right) \neq \left(\frac{2q}{s}\right)$, a contradiction. Hence, $t_{1} \neq 2q$. Consequently, $[a_{2}] \in \{2r, 2s\}$. \\
Now, let $u_{3} = x_{3} + y_{3} \sqrt{pqrs}$ be the fundamental unit of $K_{3} = \mathbb{Q}(\sqrt{pqrs})$. Then $[a_{3}] = N(u_{3}+1)$. Given $q \equiv 3 \pmod 4$, $N(u_{3}) = x_{3}^{2} - pqrs y_{3}^{2} = 1 \implies N(u_{3}+1) = 2(x_{3}+1)$. Now $pqrs \equiv 1 \pmod 4$ implies that $x_{3}$ must always be odd. This, in turn, shows that $[a_{3}] = \frac{x_{3}+1}{2}$. Denoting $\frac{x_{3}+1}{2}$ and $\frac{x_{3}-1}{2}$ by $t_{1}$ and $t_{2}$ respectively, we now try to solve $$t_{1}a^{2} - t_{2}b^{2} = 1, t_{1}t_{2} = pqrs, ab = \frac{y_{3}}{2}.$$
We first note that $\left(\frac{pr}{s}\right) = \left(\frac{ps}{r}\right)$ implies $\left(\frac{p}{r}\right) \neq \left(\frac{p}{s}\right)$. This immediate leads to $t_{1} \not \in \{p, qrs\}$. Similarly, $\left(\frac{q}{r}\right) = \left(\frac{r}{s}\right)$ implies $t_{1} \not \in \{r, qs, pqs\}$, and $\left(\frac{r}{s}\right) = \left(\frac{s}{q}\right)$ implies $t_{1} \not \in \{s, qr, pqr\}$. Finally, $\left(\frac{pr}{s}\right) = -1$ implies $t_{1} \neq pr$ while $\left(\frac{ps}{r}\right) = -1$ implies $t_{1} \neq ps$. Together, this implies $[a_{3}] \in \{q, pq, rs, prs\}$. Hence, $$\langle{\Delta_{1}, \Delta_{2}, \Delta_{3}, [a_{1}], [a_{2}], [a_{3}]\rangle} = \langle{2,p, q, r, s\rangle}.$$ Hence by Proposition \ref{prop-zan} and Lemma \ref{lem-zan2} it follows that $\mathbb{Q}(\sqrt{2p},\sqrt{2qrs})$ is a P\'{o}lya field.

\section{Compositum of Shank's cubic and Lehmer quintic : A discussion on its P\'{o}lya group and monogenicity}\label{monogenicwala}

For an integer $n \geq 1$, Shank's family of simplest cubic fields $K_{S}$ is defined by the splitting field of the irreducible polynomial $f_{n}(X) = X^{3} + (n + 3)X^{2} + nX - 1$ over the field of rational numbers $\mathbb{Q}$. The discriminant of $K_{S}$ is $(n^{2} + 3n + 9)^{2}$, whenever $n^{2} + 3n + 9$ is square-free. In a similar fashion, Lehmer's family of quintic fields $K_{L}$ is defined by the splitting field of the irreducible polynomial $g_{n}(X) = X^{5} + n^{2}X^{4} - a_{n}X^{3} + b_{n}X^{2} + c_{n}X + 1$, where $a_{n} = 2n^{3} + 6n^{2} + 10n + 10$, $b_{n} = n^{4} + 5n^{3} + 11n^{2} + 15n + 5$ and $c_{n} = n^{3} + 4n^{2} + 10n + 10$. The discriminant of $g_{n}(X)$ is given by $(n^{3} + 5n^{2} + 10n + 7)^{2}(n^{4} + 5n^{3} + 15n^{2} + 25n + 25)^{4}$.  

\smallskip

In \cite{isl}, we studied large P\'{o}lya groups in the family of Shank's simplest cubic fields and in \cite{nimish}, Lehmer's quintic fields with large P\'{o}lya groups have been investigated extensively. In this section, we briefly discuss certain possibilities of the compositum field $K_{S}K_{L}$ to have large P\'{o}lya groups and the index $1$. 

\smallskip

Let $p$ be a prime number and let $L_{p} = (p^{3} + 5p^{2} + 10p + 7)^{2}(p^{4} + 5p^{3} + 15p^{2} + 25p + 25)^{4}$ and $S_{p} = (p^{2} + 3p + 9)^{2}$. We first explore under what conditions $S_{p}$ and the first factor of $L_{p}$ are relatively prime. For that, let $q$ be a common prime factor of them. Then $q \mid p(p^{2} + 3p + 9) = p^{3} + 3p^{2} + 9p$ which implies that $q \mid [(p^{3} + 5p^{2} + 10p + 7) - (p^{3} + 3p^{2} + 9p)] = 2p^{2} + p + 7$. From this, it follows that $q \mid [2(p^{2} + 3p + 9) - (2p^{2} + p + 7)] = 5p + 11$. Again, we have $q \mid [5(p^{2} + 3p + 9) - p(5p + 11)] = 4p + 45$. Consequently, we have $q \mid [5(4p + 45) - 4(5p + 11)] = 181$ and thus $q = 181$ .

\smallskip

Similarly, if $\ell$ is a common prime divisor of both $p^{2} + 3p + 9$ and the second factor $p^{4} + 5p^{3} + 15p^{2} + 25p + 25$ of the discriminant of $g_{p}(X)$, by arguing as above, we conclude that $\ell = 541$. From this, it follows that for prime numbers $p$ satisfying the conditions in the proof of Theorem 1.2 in \cite{isl} together with the simultaneous congruence conditions $p \equiv 1 \pmod {181}$ and $p \equiv 1 \pmod {541}$, the discriminants $d_{K_{S}}$ and $d_{K_{L}}$ of $K_{S}$ and $K_{L}$, respectively, are relatively prime. Hence the discriminant $d_{K_{S}K_{L}}$ of the compositum field $K_{S}K_{L}$ equals $d_{K_{S}}^{5}d_{K_{L}}^{3}$. Moreover, since $\gcd([K_{S} : \mathbb{Q}],[K_{L} : \mathbb{Q}]) = 1$, the P\'{o}lya group $Po(K_{S}K_{L})$ of $K_{S}K_{L}$ is the direct sum of $Po(K_{S})$ and $Po(K_{L})$. Since there are infinitely many simplest cubic fields with large P\'{o}lya groups (cf. \cite[Theorem 1.2]{isl}), we conclude that for infinitely many choices of prime numbers $p$, the compositum of Shank's simplest cubic field and Lehmer's quintic field obtained from the splitting fields of $f_{p}(X)$ and $g_{p}(X)$, respectively, has large P\'{o}lya group. 

\smallskip

An interesting question at this point might be to ask whether the compositum field $K_{S}K_{L}$ is monogenic or not. In other words, whether the ring of integers $\mathcal{O}_{K_{S}K_{L}}$ is equal to $\mathbb{Z}[\alpha]$ for some $\alpha$ or not. To measure the extent of failure of a number field $F$ from being monogenic, the index of $F$ is defined as $$I(F) = \gcd\{I(\alpha) : \alpha \in \mathcal{O}_{F} \mbox{ and } F = \mathbb{Q}(\alpha)\},$$ where $I(\alpha) = [\mathcal{O}_{F} : \mathbb{Z}(\alpha)]$. Zylinski's theorem (\cite{zylin}, \cite[Proposition 4.4]{nimish}) asserts that if $[F : \mathbb{Q}] = d$ and $\ell$ is a prime divisor of $I(F)$, then $\ell < d$. In our case, we can throw more congruence relations to ensure that neither $S_{p}$ nor $L_{p}$ is divisible by any prime number smaller than $[K_{S}K_{L} : \mathbb{Q}] = 15$ for infinitely many choices of $p$. Thus we have the following theorem.
\begin{theorem}
There exist infinitely many choices of Shank's cubic fields and Lehmer's quintic fields such that their compositum is of index $1$ and has a large P\'{o}lya group. 
\end{theorem}

\begin{ack}
The authors thank their respective institutions for providing excellent facilities to carry out this work. The first author's research is supported by  CSIR (File no: 09/1026(0036)/2020-EMR-I). The discussion about the material included in Section \ref{monogenicwala} took place between the second and the last author at the restaurant Mocambo in Park Street, Kolkata over a delicious dish. We gratefully acknowledge Mocambo for playing a role behind this.
\end{ack}

\end{document}